\documentclass[11pt]{article}
\usepackage[utf8]{inputenc}
\usepackage{amsfonts,epsfig}
\usepackage[hyphens]{url}
\usepackage{breakurl}
\usepackage{comment}
\usepackage{color}
\usepackage{bbm}
\usepackage{amsmath}
\usepackage{wrapfig}
\graphicspath{{fig/}}
\usepackage{tikz}
\usepackage{relsize}
\usepackage{fancyvrb}
\usepackage{amssymb,amsmath}
\usepackage{geometry}
\usepackage{sectsty}
\usepackage{caption}
\usepackage{subcaption}
\usepackage[normalem]{ulem}
\sectionfont{\sffamily\bfseries\upshape\large}
\subsectionfont{\sffamily\bfseries\upshape\normalsize}
\subsubsectionfont{\sffamily\mdseries\upshape\normalsize}
\makeatletter
\renewcommand\@seccntformat[1]{\csname the#1\endcsname.\quad}
\makeatother\renewcommand{\bibitem}{\vskip 2pt\par\hangindent\parindent\hskip-\parindent}

\makeatletter
\def\@maketitle{%
  \begin{center}%
  \let \footnote \thanks
    {\large \@title \par}%
    {\normalsize
      \begin{tabular}[t]{c}%
        \@author
      \end{tabular}\par}%
    {\small \@date}%
  \end{center}%
}
\makeatother

\title{\bf Holes in Bayesian Statistics\footnote{For \emph{Journal of Physics G}.  We thank the participants of the Bayesian Inference in Subatomic Physics workshop and two anonymous reviewers for helpful discussion and the National Science Foundation, Office of Naval Research, and Institute for Education Sciences for partial support of this work.}\vspace{.1in}}
\author{Andrew Gelman\footnote{Department of Statistics and Department of Political Science, Columbia University.}
\and Yuling Yao \footnote{Department of Statistics, Columbia University.}
\vspace{.1in}
}
\date{13 Oct 2020\vspace{-.1in}}
\begin{document}\sloppy
\maketitle
\thispagestyle{empty}

\begin{abstract}
Every philosophy has holes, and it is the responsibility of proponents of a philosophy to point out these problems. Here are a few holes in Bayesian data analysis: (1) the usual rules of conditional probability fail in the quantum realm, (2) flat or weak priors lead to terrible inferences about things we care about, (3) subjective priors are incoherent, (4) Bayesian decision picks the wrong model,  (5) Bayes factors fail in the presence of flat or weak priors, (6) for Cantorian reasons we need to check our models, but this destroys the coherence of Bayesian inference.

Some of the problems of Bayesian statistics arise from people trying to do things they shouldn’t be trying to do, but other holes are not so easily patched. In particular, it may be a good idea to avoid flat, weak, or conventional priors, but such advice, if followed, would go against the vast majority of Bayesian practice and requires us to confront the fundamental incoherence of Bayesian inference.

This does not mean that we think Bayesian inference is a bad idea, but it does mean that there is a tension between Bayesian logic and Bayesian workflow which we believe can only be resolved by considering Bayesian logic as a tool, a way of revealing inevitable misfits and incoherences in our model assumptions, rather than as an end in itself.
\end{abstract}

\section{Overview}
Bayesian inference is logically coherent but only conditional on the assumed probability model. In the present paper we discuss several holes in the Bayesian philosophy, along with their practical implications and potential resolutions.

All statistical methods have holes and, in general, we can understand methods by understanding the scenarios where they fail (Lakatos, 1963-4).  We focus on Bayesian inference because this is the approach we use for much of our applied work and so we have an interest in deepening our understanding of it.  In addition, to the extent that coherence is a selling point of Bayesian inference, we should be aware of its limitations.

\section{Latent variables don't always exist (quantum physics)}\label{quantum}

The standard paradigm of Bayesian statistics is to define a joint distribution of all observed and unobserved quantities, with the latter including fixed parameters such as physical constants, varying parameters such as masses of different objects, and latent variables for difficult-to-measure quantities such as concentrations of a toxin within a person's liver, or, in social science, variables such as personality characteristics or political ideology that can only be measured indirectly.  The implicit assumption here is that, some joint distribution exists from which conditional and marginal probabilities can be derived using the rules of probability theory.

In quantum physics, though, there is no such joint distribution.  Or, to put it another way, quantum superposition implies that leaving a particular observation unmeasured is not the same as treating it as uncertain in a joint distribution.

Consider the canonical experiment of a light source, two slits, and a screen, with $y$ being the position where the photon hits the screen. For simplicity, think of the screen as one-dimensional, so that $y$ is a continuous random variable.

Consider four possible experiments:
\begin{enumerate}
\item Slit 1 is open, slit 2 is closed. Shine light through the slit and observe where the screen lights up. Or shoot photons through one at a time, it doesn’t matter. Either way you get a distribution, which we can call $p_1(y)$.
\item Slit 1 is closed, slit 2 is open. Same thing. Now we get $p_2(y)$.
\item Both slits are open. Now we get $p_3(y)$.
\item Now run experiment 3 with detectors at the slits. You’ll find out which slit each photon goes through. Call the slit $x$. So $x$ is a discrete random variable taking on two possible values, 1 or 2. Assuming the experiment has been set up symmetrically, you'll find that $\mbox{Pr}(x\!=\!1) = \mbox{Pr}(x\!=\!2) = 0.5$.
You can also record $y$; label its distribution $p_4(y)$.
\end{enumerate}
The difficulty arises because in experiment 4 you can also observe the conditional distributions, $p_4(y|x\!=\!1)$ and $p_4(y|x\!=\!2)$. You'll find that $p_4(y|x\!=\!1) = p_1(y)$ and $p_4(y|x\!=\!2) = p_2(y)$. You'll also find that $p_4(y) = 0.5 \,p_1(y) + 0.5\, p_2(y)$. So far, so good.

The problem is that $p_4$, which is a simple mixture distribution, is not the same as $p_3$, which has nodes and which cannot be expressed as any linear combination of $p_1$ and $p_2$.  This is a manifestation of  Heisenberg’s uncertainty principle: putting detectors at the slits changes the distribution of the hits on the screen.

This violates the rule of Bayesian statistics, by which a probability distribution is updated by conditioning on new information.  In this case, there is no joint distribution $p(x,y)$ whose marginal distribution for $y$ is $p_3$ and whose conditional distributions given $x=1$ or 2 are $p_1$ and $p_2$.  To put it another way, if we start with $p(y)=p_3(y)$, we {\em cannot} perform the marginal-conditional decomposition, $p(y)=\mbox{Pr}(x\!=\!1) p(y|x\!=\!1) + \mbox{Pr}(x\!=\!2) p(y|x\!=\!2)$.

At this point, we can rescue probability theory, and Bayesian inference, by including the measurement step in the conditioning.   Define a random variable $z$ carrying the information of which slots are open and closed, and the positions of any detectors, and then we can assign a joint distribution on $(z,y)$ and perform Bayesian inference:  Seeing the observed positions of photons on the screen gives information about the experimental and measurement process, $z$.

In essence,  the problem with applying probability theory to the slit $x$ and screen position $y$ arises because it is physically inappropriate to consider $x$ as a latent variable. The photon does not go through one slit or the other; it goes through both slits, and quantum superposition is not the same thing as probabilistic averaging over uncertainty. Mathematically, quantum mechanics can be perfectly described by Bayesian probability theory on Hilbert space given the appropriate definition of the measurable space (Hohenberg, 2010).    As long as we do not confuse additivity of elements in the vector space (i.e., $|1\rangle + |2\rangle$ is also a feasible outcome as long as $|1\rangle$ and  $|2\rangle$ are) with additivity of disjoint events under
 probability measures (i.e., $\mbox{Pr}(x=1 \text{ or }  x=2) = \mbox{Pr}(x=1)+ \mbox{Pr}(x=2)$), the quantum phenomenon is 
 simply $\{x=1 \text{ or }  x=2\} \neq \{x=|1\rangle + |2\rangle \}$, which is not surprising either probabilistically or physically.    
 
 We can redefine a simplex parameter $\theta= (\theta_1,  \theta_2)$ such that $\sum_k \theta_k=1, 0\leq \theta_k  \leq 1$ and model the superposition explicitly in the likelihood $p(y| \theta_1, \theta_2  ) =  ||\sum_{k=1}^2 \sqrt{ \theta_k}  \phi_k(y)  ||^2, $ where   $\phi_k(y) = \sqrt{(p_k(y))} \exp(i\tau_k)$ is the complex-valued wave function. Eventually we still make  inference on $\theta|y$ with uncertainty induced from finite sample observations, and still linearly average over such uncertainty  for predictive distributions $\int\! p(\tilde y|\theta)p(\theta|y)d \theta$. Bayesian inference cannot resurrect a misspecified model, but it works fine to incorporate quantum mechanics within the model.

There are two difficulties here.  The first is that Bayesian statistics is always presented in terms of a joint distribution of all parameters, latent variables, and observables; but in quantum mechanics there is no joint distribution or hidden-variable interpretation.  Quantum superposition is not merely probabilistic (Bayesian) uncertainty---or, if it is, we require a more complex probability theory that allows for quantum entanglement.

The second challenge that the uncertainty principle poses for Bayesian statistics is that when we apply probability theory to analyze macroscopic data such as experimental measurements in the biological and physical sciences or surveys, experiments, and observational studies in the social sciences, we routinely assume a joint distribution and we routinely treat the act of measurement as a direct application of conditional probability. If classical probability theory needs to be generalized to apply to quantum mechanics, then it makes us wonder if it should be generalized for applications in political science, economics, psychometrics, astronomy, and so forth.

It's not clear if there are any practical uses to this idea in statistics, outside of quantum physics. For example, would it make sense to use ``two-slit-type'' models in psychometrics, to capture the idea that asking one question affects the response to others?  Would it make sense to model macroscopic phenomena in the physical, biological, and social sciences using complex amplitudes---thus incorporating ideas of phase and coherence into applied Bayesian inference?  Could nodes, wave behavior, entanglement, and other quantum phenomena manifest in observable ways in applied statistics, motivating models that go beyond classical Boltzmann probabilities?
We have no idea.

We shall return to the two-slit experiment in section \ref{Implications}; here we emphasize that we are statisticians, not physicists, and our purpose in discussing this example is to better understand the limitations of Bayesian inference, not to claim any insight into this physics problem.

\section{Weak priors yield unreasonably strong inferences (taking posterior probabilities seriously)}\label{bets}
Consider the following not uncommon scenario.  An experiment is conducted to estimate some continuous causal effect $\theta$.  The data are summarized by an unbiased estimate, $\hat{\theta}(y)$, whose sampling distribution is approximately normal with standard error $s$.  The estimate $\hat{\theta}\pm s$ is interpreted Bayesianly under a flat prior to yield $\theta|y\sim\mbox{normal}(\hat{\theta},s)$.

Now suppose the estimate happens to be one standard error from zero, that is, $\hat{\theta}=s$.  In usual practice this would be taken as weak evidence, perhaps summarized by the 95\% posterior interval $s \pm 2s$, which overlaps zero, thus indicating that, given the information at hand, the true effect could be either positive or negative.

But the posterior distribution  could be summarized in other ways.  For example, the statement $\theta|y\sim\mbox{normal}(s,s)$ implies that $\mbox{Pr}(\theta\!>\!0|y)=0.84$; that is, there is an approximately 5/6 chance that $\theta>0$, and we should be willing to bet with 5-to-1 odds on the proposition that the treatment has a positive effect.

It would seem to be imprudent to offer such strong odds based on data that are statistically indistinguishable from noise.  To put this in a frequentist framework:  if we were to repeatedly run experiments and then offer a 5-to-1 bet on the sign of an effect whenever its estimate is 1 standerd error away from 0, we could consistently lose money.  How could pure noise result in such strong bets?

To be fair, there are problems with using betting to assess subjective probability:  every bet needs a taker, and the existence of someone willing to lay money on an event supplies information about its probability---but, even if we forget about betting and just treat that 0.84 as a probability to be used in decision making, it seems too strong.

Mathematically, 0.84 is the correct probability given the statement of the problem; if we claim that the regular application of this probability would lead to bad bets, in aggregate, this is tantamount  to arguing that the distribution of scenarios over which we are averaging is different from the prior distribution used in the model.

This example is not some sort of pathological edge case.  Estimates that are in the range of one standard error away from zero happen all the time.  And uniform priors are often taken as a default, including in our own writing (Gelman et al., 2014).  If Bayesians were to regularly lay 5-to-1 bets in this setting, they'd have lost all their money many years ago.

In practice, Bayesian statisticians deal with this problem by ignoring it, for example reporting the posterior distribution or the posterior 95\% interval and implicitly treating it as noisy and beneath notice, rather than focusing on the idea that this noisy result has, according to the model, identified the sign of $\theta$ with 5/6 probability.

The problem here is with the prior, and it's not resolved by replacing the flat prior with a proper but very weak prior, as this will still lead to a posterior probability close to 5/6 that the true effect is positive.

And this is not just an issue with extremely noisy estimates.  Consider the same example but now with an estimate $\hat{\theta}=2s$, a comfortable two standard deviations from zero.  The resulting posterior probability is now 39/40 that $\theta$ is positive, but, again, such a claim would be much too strong in practice.
By ``too strong in practice,'' we mean that this 39/40 probability would not come close to holding up, if we average over a reasonable prior distribution on true effect sizes.

Studies are commonly designed to have just enough power so that effects can just about be estimated from data.  Hence effect sizes are of the same order of magnitude as standard errors, and something like a $\mbox{normal}(0,s)$ prior may be more reasonable choice than a uniform or very weak prior on $\theta$. Studies are often optimistic, so that we'd expect true effect sizes less than $s$ much of the time.

One might say this is not a ``hole in Bayesian statistics''; it is just a implementation error to use priors that are too vague.  Sure, but then at the very least this is a mistake not just in common practice but in the recommendations in just about every textbook.  Indeed, Bayesian statisticians are often wary of strong priors because they want to ``let the data speak.''  The trouble is that letting the data speak in this way will result in many erroneously strong statements.  Paradoxically, weak priors imply inappropriately strong conclusions in certain dimensions of the posterior, a point which we have pointed out before in more complex multivariate examples; see section 3 of Gelman (1996).  Again, our point is not that a high-variance prior for a parameter of interest is necessarily wrong; rather we are emphasizing that using such a prior has implications for the applicability of the model, as it can result in strong posterior claims from what might seem like weak evidence.

\section{The incoherence of subjective priors (Bayesian updating)}\label{subjective}

Bayesians have sometimes taken the position that their models should not be externally questioned because they represent subjective belief.  The problem with this view is that, if you could in general express your knowledge in a subjective prior, you wouldn't need formal Bayesian statistics at all: you could just look at your data and write your subjective posterior distribution.

To put it another way, the inferential procedure of Bayesian statistics is to assume a prior distribution and a probability model for data and then use probability theory to determine the posterior.  But if these steps, or something approximating them, are necessary, if you {\em can't} just look at your data and come up with a subjective posterior distribution, then how is it reasonable to suppose that you could come up with an unassailable subjective distribution {\em before} seeing the data? 

We're not saying that Bayesian methods are valueless.  One virtue of the enforced consistency of Bayesian inference is that it can go in both directions.  Start with a prior and a data model, get posterior inferences, and if these inferences don't make sense, this implies they violate some aspect of your prior understanding, and you can go back and see what went wrong in your model of the world.  This is the attitude of Jaynes (2003), who takes it as a strength, not a weakness, of Bayesian inference that it makes strong claims, as these allow a model to more readily be rejected and then improved.  As O'Hagan (1994) puts it, by breaking the inference problem into smaller pieces, Bayesian inference provides a path to decision making in problems that would be too complicated to  combine data and base rate information intuitively.

Our problem is with a solipsistic doctrine of subjective priors that is sometimes expressed, in which one’s prior is considered irrefutable because it is said by definition to represent one’s subjective information. We are much more comfortable thinking about the prior (and, for that matter, the likelihood) as the product of assumptions rather than as a subjective belief.

This has implications for Bayesian practice (openness to checking and revising models) and for Bayesian theory as well:  the modeling process cannot be viewed as a closed system. The previous section demonstrated a serious, even potentially fatal, problem with noninformative or weak priors; the present section demonstrates that we cannot, except in some special cases, assume that subjectivity implies logical coherence.

\section{Bayesian decision theory can be tricky (model evaluation)}\label{decision} 
In a  decision theory based prediction framework, we equip the observable outcome $y$ and   any probabilistic prediction $p$ with a scoring rule, for example, the logarithmic scoring rule $\log p(y)$, and averaging over the actual distribution of $y$ yields the expected log predictive densities $\mathrm{E}_{y} \log p(y).$  The logarithmic scoring rule is  strictly proper, which ensures such utility, as a function of any probabilistic predictions $p$,   only attains the maximum at the true data generating process of $y$.  

Now think about experimental data in which the outcome is discrete, such as chemical measurements based on visual readings of a  test strip at discrete levels. Say we are measuring  ionizing radiation using a sloppy test kit that can only return discrete measurements: $y=0,1, \dots, 10$ mR/hr.  We also collect other predictors $x$, such as the temperature of the reaction, and the goal is to discover the relation between $x$ and $y$.  Suppose the true expected pattern follows an ordered logistic regression, $\mathrm{Pr}(y \leq k  ) =\mathrm{logit}^{-1}( c_k + \beta x) , ~  k=0,1,\dots,9$, but the practitioner chooses to fit a linear regression $y\sim \mbox{normal} (\beta x + \beta_0 , \,\sigma)$ instead. We know this model is wrong, but it could  have a higher log predictive density when fit to data.  To see this clearly,  assume that now  the test kit can  only return $y=0$ and 1 with probabilities 0.9 and 0.1 independent of $x$, so that the true data generate process Bernoulli(0.1)  has the log predictive density, or negative entropy loss,   $0.9\log(0.9)+ 0.1\log(0.1)= -0.33$, whereas a wrong normal model  $y\sim \mathrm{normal}(0, 0.3)$  achieves a higher log predictive density   $0.9\log\left(  \mbox{normal}(0|0,0.3)\right) +0.1\log\left( \mbox{normal}(1|0,0.3)\right) = -0.27$.
 
In this example, the mistake is that the calculation above mixes the continuous probability density and the discrete probability mass function, and more generally, we assign the outcome $y$  a discrete space $\{0, 1\}$ with  counting measure and simultaneously a continuous real line with Lebesgue measure. Indeed, changing the unit from mR/hr to R/hr would further inflate the log predictive density  of the wrong normal model by $\log(10^6)$, while the true discrete model remains unchanged and is always bounded above from $\log(1)=0$.

However, the measurement itself from a black box experiment can not automatically tell which space it lives in.  With observations   $y_1=1$, $y_2=0$, it is not always clear to distinguish among the following possibilities:  (1) the outcome is  continuous, and it is the test kit that induce discretization errors; (2)  the outcome has certain quantum property and is itself discrete; (3) the test kit is accurate, but it by chance runs into a few integer values in the first few trials; or (4) the test kit is nearly  accurate with certain rounding precision.
  These four schemes imply different sample spaces of the outcome. Though $\{0\}, \{1\}$ as two elements in space $\{0,1\}$,  and as elements on real line have the same notation,  they are different objects. By adopting another parameter $z$ reflecting our knowledge on what the sampling space is, there two sets of elements become
 $\{0, 1 | z=\mathrm{discrete} \}$, $\{0, 1 | z=\mathrm{continuous} \}$, akin to how conditioning on the number of open slits changes the meaning of the outcome space in the two slit experiment.  A fully Bayesian model evaluation would also average over $z$ with respect to its posterior inference when computing the utility. 

\section{Failure of Bayes factors (strong dependence on irrelevant aspects of the model)}\label{discrete}

Bayesian inference is conditional on a model.  Thus, instead of $p(\theta|y)\propto p(\theta)p(y|\theta)$, we could more formally write,  $p(\theta|y,M)\propto p(\theta|M)p(y|\theta,M)$, where $M$ is the assumed model.  But then we can go up one level of abstraction and consider a distribution $p(M)$ of possible models and then perform joint inference on models and parameters:  $p(\theta,M|y)\propto p(M)p(\theta|M)p(y|\theta,M)$.

But this joint distribution of parameters and model can be difficult to interpret, as different models can have different parameterizations so that there is no common parameter vector $\theta$ that can be usefully defined.  It  can make sense to integrate over parameters to get the marginal posterior probability of each model:  $p(M|y)=\int\! p(\theta,M|y)d\theta = p(M)\int \!p(\theta|M)p(y|\theta,M)d\theta$.  This last factor, $\int\! p(\theta|M)p(y|\theta,M)d\theta = p(M|y)$, is called the marginal likelihood of $M$, and when used to compare two models, $p(M_1|y)/p(M_2|y)$, is called the Bayes factor.

The problem with the Bayes factor is the strong dependence of the integral, $\int\! p(\theta|M)p(y|\theta,M)d\theta$, on aspects of the prior distribution, $p(\theta|M)$, that have minimal impact on the conditional posterior distribution, $p(\theta|y,M)$.   From a practical standpoint, weak or noninformative priors can often be justified on the grounds that, with strong data, inferences for quantities of interest are not seriously affected by the details of the prior---but this is not the case for the marginal likelihood.

As with all the other problems discussed in this article, there is nothing wrong with the Bayes factor from a mathematical perspective.  The problem, when it comes, arises from the application.

There are three sorts of examples where Bayes factors can be applied, two where the approach can make sense (if performed carefully) and one where it doesn't.

{\em Example 1:  Continuous parameter.}  The simplest example of Bayes factors is when the model index $M$ represents a continuous parameter.  For example, consider a logistic regression with two coefficients, $\mbox{Pr}(y_i=1)=\mbox{logit}^{-1}(a + bx_i)$, for $i=1,\dots,n$, with the model implicitly conditional on the values of the predictor $x$.  Suppose we identify the intercept, $a$, as the continuous model index $M$, and the slope $b$ as the parameter $\theta$.  The we can work with the joint distribution $p(\theta,M|y)$, the marginal model probability distribution $p(M|y)$, the marginal likelihood $p(y|M)$, the model-averaged posterior $p(\theta|y)$, and so forth.  There is no problem at all; it's just Bayesian inference.

{\em Example 2:  Truly discrete model.}  Another setting where Bayes factors can make sense is in a discrete mixture model.  For example, $M=1$ or 2, depending on whether a person was or was not at the scene of a particular crime, with data $y$ representing some sort of surveillance measurement characterized by some parameters $\theta$.  Here, if we want to compute the posterior probability of the person being at the crime, we need to integrate over $\theta$.

{\em Example 3:  Model choice.}  Bayes factors run into trouble when used to compare discrete probabilities of continuous models.  For example, consider the problem of comparing two logistic regressions, one with just an intercept ($M_1: \ \mbox{Pr}(y_i=1)=\mbox{logit}^{-1}(a)$) and the other including a slope as well ($M_2: \ \mbox{Pr}(y_i=1)=\mbox{logit}^{-1}(a + bx_i)$).  The Bayes factor comparing the two models will depend crucially on the prior distribution for the slope, $b$, in the second model.  For example, suppose the models are on unit scale so that we expect the parameters $a$ and $b$ to be of order 1.  Then switching from the priors $b|M_2\sim\mbox{normal}(0,10)$ to $b|M_2\sim\mbox{normal}(0,1000)$ will multiply the Bayes factor by approximately 100, even though this change will have essentially no influence on the parameters $a$ and $b$ within the model.  In real-world physics modeling, this sort of indeterminacy should not occur, as our understanding of the units of the problem should allow us to set a reasonable prior to within an order of magnitude, but for models with many parameters, smaller factors can still accumulate to the extent that the Bayes factor can still depend strongly on arbitrary settings of the priors.
In the literature on the Bayes factor, solutions to this problem have been proposed, along the lines of setting some conventional prior on parameters within each of the candidate models.  But we have never found these arguments convincing---and, to the extent that the are convincing, they run counter to the goal of coherence, by which prior distributions are intended to represent actual populations to be averaged over when evaluating inferences.

\section{Cantor's corner (model updating)}\label{cantor}
A model does what it does, then we look for problems; having done that, we expand the model and push it harder, until it breaks down and we need to replace it.  This is basic philosophy of science since Lakatos (1968), and it maps directly to Cantor's diagonal argument, a proof that the real numbers are uncountable.

Here's the story in Ascii art:

\begin{small}\begin{verbatim}
     Model 1 . . . X
     Model 2 . . . . . . X
     Model 3 . . . . . . . . X
     Model 4 . . . . . . . . . . . . . . . . X
     Model 5 . . . . . . . . . . . . . . . . . . . X
     . . . 
\end{verbatim}\end{small}

\noindent
For each model, the dots represent successful applications, and the X represents a failure, a place where the model needs to be improved. The X's form Cantor’s diagonal, and the most recent X is Cantor’s corner. Cantor's diagonal argument, taken metaphorically, is that we can never fully anticipate what is needed in future models. Or, to put it another way, the sequence of problems associated with the sequence of X's cannot, by definition, be fit by any single model in the list.

Scientific research is all about discovery of the unexpected: to do research, you need to be open to new possibilities, to design experiments to force anomalies, and to learn from them. The sweet spot for any researcher is at Cantor's corner.  As such, direct Bayesian inference---either the discrete version with model averaging criticized in section \ref{discrete}, or a more desirable approach of continuous model averaging---is, ultimately, impossible in that it would require the anticipation of future steps along the diagonal, which by construction cannot be done. 
There is a connection here to the classical idea of sieves:  models that expand with the data (Grenander, 1981, Geman and Hwang, 1982).

One way to fix ideas here is to imagine a computer program or artificial intelligence that would perform statistical inference on arbitrary datasets.  One could stock this program with some classes of useful models such as linear regressions, logistic regressions, splines, tree models, and so forth, along with some rules that would allow us to put these pieces together into arbitrarily---countably many---larger models, for example using interactions, hierarchical modeling, network connections, deep nets, nets of nets, etc.  The computer program could then try out a series of models on any given dataset, using cross validation and more complicated evaluation and combination rules such as random forests.  Each of these steps allows a model, or class of models, to last longer, to work deeper in the class of applied problems, but at some point the model would suffer serious misfit and stop performing better.  To be an effective general-purpose statistics problem solver, the AI would need to include a module for model criticism:  a sort of ``red team'' to detect problems which at some point would require ripping out many stitches of modeling before future progress can be made.  This idea is similar to generative adversarial networks in machine learning (Goodfellow et al., 2014); our point is that the potentially unlimited nature of Cantorian updating reveals another way in which Bayesian inference cannot hope to be coherent in the sense of being expressible as conditional statements from a single joint distribution.

\section{Implications for statistical practice}\label{Implications}
We continue by discussing how we can incorporate into Bayesian workflow our understanding of the holes in Bayesian inference.

First consider the problem in section \ref{quantum} of the two-slit experiment, again emphasizing that our purpose here is not to contribute to understanding of this physics problem but rather to gain insight into the challenges of statistical inference in problems of quantum indeterminacy. Bayesian (Boltzmann) probability works in experiments 1, 2, and 4, but not in experiment 3, where there is no joint distribution for the slit $x$ and the screen position $y$.  One way to rescue this situation and make it amenable to standard Bayesian analysis is to fully respect the uncertainty principle and only apply probability models to what is measured in the experiment. Thus, if both slits are open and there are no detectors at the slits, we do not express $x$ as a random variable at all. In a predictive paradigm, inference should only be on observable quantities (Bernardo and Smith, 1994). This is related to the argument of Dawid (2000) that it can be mistake to model the joint distribution of quantities such as potential outcomes that can never be jointly observed.

Bayes' rule by default aggregates uncertainty of potential outcomes under different measurements through a linear mixture:   the marginal distribution of $y$ is a mixture of the conditional distributions under the prior distribution: $p(y)= \int\! p(y|x)p(x) dx$, and  the posterior predictive distribution of  future outcome $\tilde{y}$ is a mixture under the posterior distribution of $x$: $p(\tilde y | y)= \int\! p(\tilde y| x) p(x| y ) d x.$ However, the linear mixture is just one operator on the space of predictive distributions, among many others such as convolution, multiplication of densities, or superposition. That said, restricting the average form to a linear mixture is not unique to Bayesian inference: a confidence interval from Neyman–Pearson theory also makes a claim that is linearly averaged over all possible scenarios under consideration.

We explained in section  \ref{quantum}  that the Bayesian inference can be rescued by careful modeling.  A challenge of this model-the-observables approach is that it requires some external knowledge. In particular, the mathematics of wave mechanics are required to come up with the distribution $p_3(y)$ as a superposition of $p_1(y)$ and $p_2(y)$, with its counterintuitive interference pattern.  In addition, there is the awkwardness that $x$ can be measured for some conditions of the number of open slits but not others. 

But, from the standpoint of usual statistical practice, it is not at all unusual for different models for the same outcome to have parameters that are not directly comparable.  For example, in pharmacology, one might model the concentration of a compound in the bloodstream (an observable outcome) jointly with its unobservable concentrations within internal compartments of the body, or simply fit a curve directly to the blood measurements.  Similarly, an economist might fit a ``structural'' model including latent variables representing individual utilities, or a ``reduced-form'' or phenomenological model to data on purchases.  In the pharmacology example, the latent variables seem uncontroversial and we could hope that the direct model of blood concentration could be expressed mathematically as an integral over the unknown internal concentrations, thus following Bayes' rule.  But in the economics example, it has often been argued that preferences do not exist until they are measured (Lichtenstein and Slovic, 2006), hence any joint model would depend on the measurement protocol.

More generally, we can think about the question of {\em what} to model as part of the modeling process.  Without an explicit model, it is the default linear mixture  ($0.5\, p(y|x=1) +0.5 \,p(y|x=2)$) that misfits the quantum observations.  Likewise, without an explicit prior, the default uniform prior on all parameters can yield undesired results.

In section \ref{bets} we considered problems of noninformative or weakly informative priors, but arguably the key mistake in such settings comes from the decision not to include relevant prior information on a parameter $\theta$, a decision which in the Bayesian setting corresponds to a refusal to consider $\theta$ as a draw from some real or fictive population of possible parameters. In the example of section \ref{bets}, these would be the possibility that individual treatment effects vary across sub population, whereas the  posterior inference $\theta|y$ ignores such unmodeled epistemic uncertainty and has to converge to the average treatment effect regardless,  thereby resulting in an overconfident bet for individual outcomes.   By contrast, point estimation procedures such as maximum likelihood do not even try to this end.

To be more constructive, in any given applied problem we can consider potential decision rules such as betting on the posterior probability that $\theta>0$ and, if its implications under hypothetical repeated sampling seem undesirable, we can consider this as a prior predictive check (Box, 1980) or ``device of imaginary results'' (Good, 1950) that reveals additional beliefs that we have which are inconsistent with our assumed model.  This revelation of incoherence does not immediately tell us what to do, but it motivates a more careful engagement with our assumptions. In particular, the model itself does not know how it will be used in the decision problem: A time series model good at one-day ahead prediction does not necessarily yield the optimal one-year ahead prediction.     

To put it another way, it is a fundamental principle of Bayesian inference that statistical procedures do not apply universally; rather, they are optimal only when averaging over the prior distribution.  This implies a proof-by-contradiction sort of logic, relating to the discussion of subjective priors in section \ref{subjective} by which it should be possible to deduce properties of the prior based on understanding of the range of applicability of a method.

The example in section \ref{decision} makes it clear that decision theory itself is not immune to subjective model building,  or the  subjective basis of knowledge.  We need prior knowledge, an interpretation of how the observations are generated, to extrapolate the measurements to the underlying constructs of interest.  In many problems it is acceptable to use a linear regression for people's height with a normal error even though we know the measurement has to be positive and bounded, and knowing  that all  measurements recorded by float numbers are in essence discrete does not prevent us applying continuous probability  densities.  The choice of sampling space is part of the model assumptions and approximations we are willing to make, and in principle we can further check how sensitive the analysis is to these assumptions as part of the model. 

What about the failure of Bayes factors, as discussed in section \ref{discrete}?  At the most immediate level of statistical practice, we can use more robust Bayesian model averaging techniques such as stacking (Yao et al., 2018) so as to use the fit of models to data to obtain better predictive performance without the pathologies of the marginal likelihood.  At a more fundamental level, we are again seeing the consequences of a refusal to model a parameter $\theta$, which has unfortunate implications when integrating over $\theta$ in the prior distribution.

It would be awkward, however, to simply insist on realistic priors for all parameters in our Bayesian inferences.  Realistically it is not possible to capture all features of reality in a statistical model.  The relevant point here is that the refusal to model some part of our problem is a choice that narrows the range of what can be done going forward.  There is no fundamental problem with Bayes factors and Bayesian model averaging, but there is a problem with integrating over a parameter that is essentially unmodeled.

Alternatively we view the drawback of Bayesian model averaging as the refusal to model how individual models are combined.
Bayesian model averaging is indeed a form of averaging, but it assumes that one of the candidate models is true, with individual models aggregated according to the probability of each model being that true model. With more and more data, the  Bayesian model averaging  will asymptotically only pick one model and becomes marginal likelihood based model \emph{selection}, which is inappropriate in a model-open environment in which none of the posited models are correct (Yao et al., 2018).
The default assumption behind the notation $p(M|y)$, that there is one true model and it is in the set under consideration, is akin to the implicit assumption that only one slit is open in the protocol induced by the notation $p(x|y)$ in the two-slit experiment. Although the Bayesian inference from finite sample induces nonzero epistemic uncertainty on $M$ or $x$, the underlying decision space is essentially discrete, either allowing individual models $M_k$ to be binarily true or false, or the slit to be open or closed.  Through additional modeling in the two-slit experiment,  we are expanding the binary decision $x$ into a two-dimensional simplex $\theta$ and likewise expanding the binary model of truth $\{M_k: M_k=0 \text{ or } 1, \sum_{k=1}^K M_k= 1\}$ into a $K$-dimensional simplex $\{w_k: 0 \leq w \leq 1, \sum_{k=1}^Kw_k= 1\}$. The enlarged decision space is a continuous extension of the original one, rendering more flexibility to express the data generating process.

Finally, as discussed in section \ref{cantor}, Cantor's diagonal argument points to the essential ongoing nature of these problems:  there will always be more aspects of reality that have not yet been included in any model.  Lest this be seen as an endorsement of some sort of flabby humanism of the ``scientists are people too'' variety, let us emphasize that the same issues arise with extrapolation and generalization from automatic machine learning algorithms.

\section{Discussion}

There are various reasons to care about holes in Bayesian statistics.  Most directly, we want to avoid these holes in our applied work, and we want to be aware of alternative approaches, going beyond default noninformative models or naive overconfidence in subjective models.  Stepping back, we should be appropriately wary of statistical inferences produced using any method, Bayesian or otherwise, as the Cantorian argument implies that coherence will never be within reach.

The holes discussed in this article are not simply the errors of inference resulting from inevitable imperfections in model misspecification.  Rather, each hole in its own way represents a potentially catastrophic failure:  in the two-slit experiment, the naive application of conditional and joint probability leads to a mixture distribution rather than the actual result with nodes; in the flat-priors example, the objectionably strong claim that $\Pr(\theta>0)=0.84$ can arise from noisy data; the incoherence of Bayesian updating eliminates any theoretical basis for a behavioristic interpretation of the prior distribution; Bayes factors cannot be rescued, even approximately, by the use of conventional priors; and the essence of Cantor's corner is that any model, if used in a changing world, will eventually need radical overhaul.  These holes represent real challenges to any philosophy of automatic Bayesian inference.

As we have discussed elsewhere (Gelman and Robert, 2013), ``in the popular and technical press, we have noticed that `Bayesian' is sometimes used as a catchall label for rational behavior.''  But ``rationality (both in the common-sense and statistical meanings of the word) is complex. At any given time, different statistical philosophies will be useful in solving different applied problems. As Bayesian researchers, we take this not as a reason to give up in some areas but rather as a motivation to improve our methods: if a non-Bayesian method works well, we want to understand how.''   Awareness of the holes in Bayesian statistics should allow us to be better Bayesians in the short run, while pointing to research directions for the future.

There are various levels of Bayesian statistical practice.   \emph{Bayesian inference} goes from likelihood and prior to posterior,   a deterministic procedure followed by  Bayes' theorem after one model is set. \emph{Bayesian logic} contains both Bayesian inference and decision theory.    \emph{Bayesian data analysis}  involves  three steps of model building, inference, and model checking and improvement.  A complete \emph{Bayesian workflow} refers to Bayesian data analysis for a sequence of models,  not just for the purpose of model choice or model averaging but more importantly to understand these models.

When discussing logical gaps, our target of criticism is not Bayesian inference as a prior-to-posterior mapping or Bayesian inference as model-based reasoning; rather, we are pointing out flaws in a static form of Bayesian inference in which the model is a fixed, passive receptacle for the data.  Here is the challenge:  It is fine to say that we should be flexible and learn from our models' failures, and indeed that is how we view Bayesian workflow (Gabry et al., 2019), but this directly conflicts with one of the usual justifications for Bayesian inference, which is its logical consistency (Cox, 1946, Savage, 1954).  Again, the recognition of these holes should not be taken as a reason to abandon Bayesian methods or Bayesian reasoning, but rather as a motivation to better understand then to improve modeling and inferential workflow, and to better understand the dependence of inferences on model assumptions.

 Bayesian statistics has no holes in theory, only in practice.    When the model and all model assumptions are taken for granted, the Bayesian procedure is unique.  But we prefer to break such enforced coherence by questioning the validity of the model and  assumptions.  
 On the other hand, perhaps these practical holes imply the existence of a theoretical hole, if theoretical statistics indeed is the theory of applied statistics.
\section*{References}

\noindent

\bibitem Bernardo, J. M., and  Smith, A. F. M. (1994). {\em Bayesian Theory}. New York: Wiley.

\bibitem Box, G. E. P. (1980).  Sampling and Bayes inference in scientific modelling and robustness.  {\em Journal of the Royal Statistical Society A} {\bf 143}, 383--430.

\bibitem Cox, R. T. (1946). Probability, frequency and reasonable expectation. {\em American Journal of Physics} {\bf 14}, 1--10.

\bibitem Dawid, A. P. (2000).  Causal inference without counterfactuals.  {\em Journal of the American Statistical Association} {\bf 95}, 407--424.

\bibitem Feller, W. (1950). {\em An Introduction to Probability Theory and Its Applications}. New York: Wiley.

\bibitem Gabry, J., Simpson, D., Vehtari, A., Betancourt, M., and Gelman, A. (2019).  Visualization in Bayesian workflow (with discussion). {\em Journal of the Royal Statistical Society A} {\bf 182}, 389--402.

\bibitem Gelman, A. (1996). Bayesian model-building by pure thought: Some principles and examples. {\em Statistica Sinica} {\bf 6}, 215--232.

\bibitem Gelman, A., Jakulin, A., Pittau, M. G., and Su, Y. S. (2008).  A weakly informative default prior distribution for logistic and other regression models. {\em Annals of Applied Statistics} {\bf 2}, 1360--1383.

\bibitem Gelman, A., and Robert, C. (2013).  ``Not only defended but also applied'': The perceived absurdity of Bayesian inference (with discussion). {\em American Statistician} {\bf 67}, 1--17.

\bibitem Geman, S., and Hwang, C. R. (1982).  Nonparametric maximum likelihood estimation by the method of sieves.  {\em Annals of Statistics} {\bf 10}, 401--414.

\bibitem Good, I. J. (1950).  {\em Probability and the Weighing of Evidence}.  New York:  Hafner.

\bibitem Goodfellow, I., Pouget-Abadie, J., Mirza, M., Xu, B., Warde-Farley, D., Ozair, S., Courville, A., and Bengio, Y. (2014). Generative adversarial nets. {\em Advances in Neural Information Processing Systems} {\bf 27}, 2672--2680.

\bibitem Grenander, U. (1981).  {\em Abstract Inference}.  New York:  Wiley.

\bibitem Griffiths, R. B. (2002). {\em Consistent Quantum Theory}.  Cambridge University Press.

\bibitem Hohenberg, P. C. (2010).  An introduction to consistent quantum theory.  {\em Reviews of Moden Physics} {\bf 82}, 2835--2844.

\bibitem Jaynes, E. T. (2003).  {\em Probability Theory:  The Logic of
Science}. Cambridge University Press.

\bibitem Lakatos, I. (1963-4).  Proofs and refutations.  {\em British Journal for the Philosophy of Science} {\bf 14}, 1--25, 120--139, 221--243, 296--342.

\bibitem Lakatos, I. (1968). Criticism and the methodology of scientific research programmes.  {\em Proceedings of the Aristotelian Society} {\bf 69}, 149--186.

\bibitem Lichtenstein, S., and Slovic, P. (2006).  {\em The Construction of Preference}.  Cambridge University Press.

\bibitem O'Hagan, A. (1994).  {\em Kendall's Advanced Theory of Statistics, Bayesian Inference}.  New York:  Wiley.

\bibitem Savage, L. J. (1954).  {\em The Foundations of Statistics}.
New York:  Dover.

\bibitem Yao, Y., Vehtari, A., Simpson, D., and Gelman, A. (2018). Using stacking to average Bayesian predictive distributions (with discussion).  {\em Bayesian Analysis} {\bf 13}, 917--1007.

\pagebreak

\section{Appendix:  Things that seem like problems with Bayesian inference but aren't}

Bayesian inference has been criticized from many directions.  Here we wish to briefly supplement the holes described in the main body of this paper, with some other criticisms that we see as less fundamental.

\subsection{Bayes as a pretty but impractical theory}
A famous probabilist once wrote that Bayesian methods ``can be defended but not applied'' (Feller, 1950).  There is a pragmatist appeal to the idea that Bayesian inference is doomed by its own consistency, but ultimately we find this particular anti-Bayesian argument unconvincing, for two reasons.  First, much as changed since 1950, and Bayesian methods can and are applied in many areas of science, engineering, and policy.  Second, the consistency of Bayesian inference is only a problem if you are inflexible.  As discussed in section \ref{bets}, if a posterior distribution yields unappealing inferences, this is an opportunity to interrogate the model and add information as appropriate.  A probability model is a tool for learning, not a suicide pact.

In making his comparison, Feller was defining Bayesian statistics by its limitations while crediting classical hypothesis testing with the 1950 equivalent of vaporware.  He perhaps leapt from the existence of a philosophical justification for Bayesianism to an assumption that philosophical arguments were the {\em only} justification for Bayesian inference.

\subsection{Bayes as an automatic inference engine}
In bygone days, Bayesian inference was criticized as too impractical to apply in any but the simplest problems.  With the advent of algorithms such as the Gibbs sampler, Metropolis, and Hamiltonian Monte Carlo, along with black-box implementations such as Bugs and Stan, we sometimes hear the reverse criticism, that Bayesian models can be fit too easily and thoughtlessly, in contrast to classical inference, where the properties of statistical models have to be fit laboriously, one at a time.  We agree that the relative ease of Bayesian model fitting puts more of a burden on model evaluation---with great power comes great responsibility---but we take this as a criticism of rigid, fixed-model Bayesian inference, not applying to the open version discussed in section \ref{cantor}.

\subsection{Bayesian inference as subjective and thus nonscientific}
We have proposed to replace the words ``objective'' and ``subjective'' in statistics discourse with broader collections of attributes, with objectivity replaced by {\em transparency}, {\em consensus}, {\em impartiality}, and {\em correspondence to observable reality}, and subjectivity replaced by {\em awareness of multiple perspectives} and {\em context dependence} (Gelman and Hennig, 2017).  With this in mind, we believe that concerns about the subjectivity of Bayesian prior distributions and likelihoods can be addressed  first, by grounding these choices in correspondence to observable reality (i.e., prior data); and, second, by making these choices transparently and with an awareness of contexts.  Thus, rather than criticizing Bayesian modeling as subjective, we think of modeling choices as embedded in a larger workflow of information gathering, synthesis, and evaluation.

\end{document}